\documentclass[11pt,reqno]{amsart}
\allowdisplaybreaks
\usepackage{graphicx}
\newtheorem{thm}{Theorem}[section]
\newtheorem{cor}[thm]{Corollary}
\newtheorem{lem}[thm]{Lemma}
\newtheorem{prop}[thm]{Proposition}
\newtheorem{rmk}[thm]{Remark}
\numberwithin{equation}{section}

\newcommand{\hs}{\hspace}

\newcommand{\td}{\tilde}

\newcommand{\fr}{\frac}
\newcommand{\ed}{\end{document}}
\newcommand{\be}{\begin{equation}}
\newcommand{\ee}{\end{equation}}

\newcommand{\lagl}{\langle}
\newcommand{\ragl}{\rangle}
\newcommand{\lmx}{\left(\begin{matrix}}
\newcommand{\rmx}{\end{matrix}\right)}
\newcommand{\ldt}{\left|\begin{matrix}}
\newcommand{\rdt}{\end{matrix}\right|}

\newcommand{\rank}{{\rm rank\,}}

\newcommand{\vfi}{\varphi}

\newcommand{\bbr}{{\mathbb R}}

\newcommand{\ba}{\begin{array}}
\newcommand{\ea}{\end{array}}
\newcommand{\bald}{\begin{aligned}}
\newcommand{\eald}{\end{aligned}}
\newcommand{\bea}{\begin{eqnarray}}
\newcommand{\eea}{\end{eqnarray}}

\newcommand{\id}{{\rm id}}

\begin{document}

\title[Isometric Immersions into $S^k \times H^{n+p-k}$ of higher codimension]{Isometric immersions of higher codimension\\
into the product ${\mathbf S^k \times H^{n+p-k}}$ of sphere and hyperboloid}
\thanks{Research supported by grants of NSFC (No. 11171091).}

\author[X. X. Li and T. Q. Zhang]{Xingxiao Li $^*$ and Tianqun Zhang}
\thanks{$^*$ Corresponding author.}
\address{Department of Mathematics, Henan Normal University, Xinxiang 453007, Henan, P.R.China}%
\email{xxl@henannu.edu.cn}%


\begin{abstract}
In this paper, we
obtain a sufficient and necessary condition for a simply connected
Riemannian manifold $(M^n, g)$ to be isometrically immersed, as a submanifold with codimension $p\geq 1$, into
the product $S^k\times H^{n+p-k}$ of sphere and hyperboloid.
\end{abstract}
\maketitle

\section{Introduction}

It is well known that the Gauss and Codazzi equations are
necessary conditions for a Riemannian manifold $(M^n, g)$ of dimension
$n$ to be locally
isometrically immersed into an arbitrary $(\bar n)$-dimensional Riemannian manifold
$(\bar M ^{\bar n},\bar g)$ as a hypersurface. Then a classical problem is: when can a Riemannian
manifold $(M^n, g)$, with prescribed first and second fundamental
forms satisfying the Gauss and Codazzi equations, be isometrically
immersed into a Riemannian manifold $(\bar M ^{\bar n},\bar g )$?
This problem was solved by Bonnet for $\bar M =\bbr^3$ (see \cite{bo} or
\cite{k-n}). For arbitrary dimensions and co-dimensions, Tenenblat gave in
1971 an elementary proof of the fundamental theorem for
immersions into the Euclidean space $\bbr^{n+p}$ (see \cite{te}). In
2005, B\"ar etc gave a proof for hypersurfaces in $\bbr^{n+1}$ which
can be directly extended to
semi-Riemannian manifolds (see \cite{b-g-m}).

When $\bar M ^{\bar n}$ is a space form, the Gauss and
Codazzi equations are intrinsicly defined on submanifolds, and they
are also sufficient for an n-dimensional simply connected Riemannian
manifold to be isometrically immersed into $\bar M ^{\bar n}$ (see
\cite{sp}). Recently, B. Daniel obtained necessary and sufficient conditions
for an $n$-dimensional simply connected Riemannian manifold to be isometrically
immersed into $S^n\times \bbr$ or $H^n\times \bbr$ in terms of its first
and second fundamental forms and the projection of the vertical
vector field (see \cite{da}). If the codimension is greater than one, then
the fundamental equations of submanifold involve an additional
equation, that is, the Ricci equation. Q. Chen and Q. Cui
obtain in \cite{c-c} sufficient and necessary conditions for a simply connected submanifold
$(M^n, g)$ to be isometrically immersed into
$S^m\times \bbr$ and $H^m\times \bbr$, respectively. In \cite{ko}, D. Kowalczyk gives sufficient and necessary conditions for a simply
connected hypersurface $(M^n, g)$ in $S^k\times H^{n-k+1}$.

 In this paper, we will extend this result to $S^k\times
  H^{n+p-k}$ with general codimensions. Contrast to Chen-Cui in \cite{c-c}, our theorem only uses
the global invariants. When $p=1$, our
treatment easily recovers the main theorem of D. Kowalczyk.

As we know, on the Riemannian product
  $S^k\times H^{n+p-k}$, there exists a natural symmetric
  $(1,1)$-tensor $\bar{\psi}$, called the product
  structure of $S^k\times H^{n+p-k}$ satisfying $\bar{\psi}^2=\id$ $(\bar{\psi} \neq \pm \id)$ and
  $\tilde{D}\bar{\psi}=0$, where $\tilde{D}$ is the induced connection of the Levi-Civita
  connection $\bar D$ on $S^k\times H^{n+p-k}$.

  Let $\varphi:(M^n,g)\rightarrow(S^k\times
  H^{n+p-k}, \bar g $) be an isometric immersion
  of codimension $p$ with the second fundamental form $h$. Denote by $TM^n$ (resp. $T^\bot_\vfi M^n$) the tangent (resp. normal)
  bundle of $\vfi$. Then we have the following bundle decomposition:
  $$
  \vfi^*T(S^k\times H^{n+p-k})=\vfi_*TM^n\oplus T^\bot_\vfi M^n.
  $$
Write
\be\label{1.1}\bar{\psi}(\vfi_*X)=\vfi_*f(X)+\bar{u}(X),\quad
\bar{\psi}(\xi)=\vfi_*(\bar{U}(\xi))+\bar{\lambda}(\xi),
\ee
for all $X\in TM^n$, $\xi\in T^{\perp}M^n$, where
    $$f:TM^n
    \rightarrow TM^n, \quad\bar{u}:TM^n\rightarrow T^{\perp}M^n,\quad\bar{U}:T^{\perp}M^n
    \rightarrow
    TM^n,\quad \bar{\lambda}:T^{\perp}M^n\rightarrow T^{\perp}M^n
    $$
are bundle homomorphisms. If $R$ (resp. $R^\bot$) is the curvature (resp. the normal curvature) of $\vfi$ then,
in the case of $\bar M^{n+p}=S^k\times H^{n+p-k}$, we can rewrite respectively the
    equations of Gauss, Codazzi and Ricci as follows (see Proposition \ref{p1}):
\begin{align*}
&R(X,Y)Z=A_{h(Y,Z)}X-A_{h(X,Z)}Y\\
&\hs{2.2cm}+\frac{1}{2}(g(f(Y),Z)X+g(Y,Z)f(X)-g(f(X),Z)Y-g(X,Z)f(Y)),\\
&2((\tilde{D}_{X}h)(Y,Z)-(\tilde{D}_{Y}h)(X,Z))=g(Y,Z)\bar{u}(X)-g(X,Z)\bar{u}(Y),\\
&R^\bot(X,Y)\xi=h(A_{\xi}Y,X)-h(A_{\xi}X,Y),
\end{align*}
where $X,Y,Z\in \Gamma(TM^n)$, $\xi\in \Gamma(T^\bot_\vfi M^n)$, $A_\xi:TM^n\to TM^n$ is the Weingarten map of $\vfi$ and $\td D$ is the induced connection.

To find the sufficient conditions for an $n$-dimensional Riemannian manifold $(M^n,g)$ to be isometrically immersed into $S^k\times H^{n+p-k}$, we need some suitably chosen metric vector bundle $E\to M^n$ in the place of $T^\bot_\vfi M^n$ and some $E$-valued symmetric $2$-form $\sigma$ in the place of $h$.

Now for any metric vector bundle $E\to M^n$ with metric $\lagl\cdot,\cdot\ragl_E$ and with a given metric connection $D^E$ and an $E$-valued symmetric $2$-form $\sigma$, that is, $\sigma\in\Gamma(T^*M^n\odot T^*M^n\otimes E)$, define the metric vector bundle sum $\bar E:=TM^n\oplus E$ and the corresponding metric on $\bar E$ is denoted by $\lagl\cdot,\cdot\ragl_{\bar E}$. Using the connections $D$ on $TM^n$, $D^E$ on $E$ and the $2$-form $\sigma$, we can introduce one unique metric connection $D^{\bar E}$ on $\bar E$ as follows:
\begin{align}
D^{\bar E}_XY=&D_XY+\sigma(X,Y),\quad X\in TM^n,\ Y\in\Gamma(TM^n),\label{1.1-1}\\
D^{\bar E}_X\xi=&-A^\sigma_\xi(X)+D^E_X\xi,\quad X\in TM^n,\ \xi\in \Gamma(E),\label{1.1-2},
\end{align}
where, for each $\xi\in \Gamma(E)$, the bundle map $A^\sigma_\xi:TM^n\to TM^n$
is defined by
\be\label{1.10}
g(A^\sigma_\xi(X),Y)=\lagl\sigma(X,Y),\xi\ragl_E,\quad \forall X,Y\in TM^n.\ee

In the present paper, we aim at proving the following main theorem.

{\thm\label{main} Let $(M^n,g)$ be a simply connect Riemannian
manifold of
     dimension $n$. If there exist a Riemannian vector bundle
     $E \rightarrow M^n$ of rank $p(\geq 1)$ with a metric connection
     $D^E$ of curvature $R^E$, an $E$-valued symmetric form
     $\sigma\in\Gamma(T^{*}M^n\odot T^{*}M^n\otimes E)$ and a symmetric
      bundle endomorphism $\psi(\neq\pm\id)$ of the metric sum bundle $\bar E=TM^n\oplus E$, with the connection $D^{\bar E}$ defined by \eqref{1.1-1} and \eqref{1.1-2}$:$
      \be
      \psi=\lmx
     f&U\\u&\lambda\rmx\in {\rm End}(TM^n\oplus E)\label{1.5}
     \ee
    where
$$
f:TM^n\rightarrow TM^n, \quad u:TM^n\rightarrow E,\quad
U:E\rightarrow TM^n,\quad \lambda:E\rightarrow E
$$
are bundle homomorphisms,  such that for all
$X,Y,Z\in \Gamma(TM^n)$, $\xi\in \Gamma(E)$
\begin{align}
&\psi^2=\id_{TM^n\oplus E},\quad \td D\psi=0,\label{1.6}\\
&R(X,Y)Z=A^\sigma_{\sigma(Y,Z)}X-A^\sigma_{\sigma(X,Z)}Y\nonumber\\
&\hs{2.2cm}+\frac{1}{2}(g(Y,Z)f(X)-g(X,Z)f(Y)+g(f(Y),Z)X-g(f(X),Z)Y),\label{1.7}\\
&2\left((\tilde{D}_{X}\sigma)(Y,Z)-(\tilde{D}_{Y}\sigma)(X,Z)\right)=g(Y,Z)u(X)-g(X,Z)u(Y),\label{1.8}\\
&R^E(X,Y)\xi=\sigma(A^\sigma_{\xi}Y,X)-\sigma(A^\sigma_{\xi}X,Y),\label{1.9}
\end{align}
where $\td D$'s are the induced connections and, for each $\xi\in\Gamma(E)$, the bundle map $A^\sigma_\xi:TM^n\to TM^n$ is defined by \eqref{1.10}.
Then there exit some $k$, $1\leq k\leq n+p-1$, and an isometric
immersion $\varphi: M^n \rightarrow S^k\times
     H^{n+p-k}$, unique up to isometries of $S^k\times
     H^{n+p-k}$, and an isometric bundle homorphism $\Phi\in {\rm Hom}(E,T^{\perp}_\vfi M^n)$ satisfying
\be\label{1.11}\bar\psi\circ\varphi_{ \ast}=\varphi_{ \ast}f+\Phi\circ u,\quad
\bar\psi\circ\Phi=\varphi_{\ast}U+\Phi\circ\lambda.\ee
Moreover, $h:=\Phi\circ\sigma$ is the second fundamental form of
$\varphi$.

Furthermore, the above conditions \eqref{1.6}--\eqref{1.9} are also necessary for the existence of an isometric immersion $\vfi:M^n\to S^k\times H^{n+p-k}$.}

\rmk\rm We would like to remark here that, by a minor modification of the same argument, one can prove that our main theorem remains true if one replaces the ambient product space $S^k\times H^{n+p-k}$ by a direct product $\bar M^k(c_1)\times \bar M^{n+p-k}(c_2)$ of any two real space forms $\bar M^k(c_1)$, $\bar M^{n+p-k}(c_2)$ with constant curvatures $c_1,c_2$, respectively.

The present paper is organized as follows: In the next section, we shall recall some basic facts that are to be used later, including the some concepts and formulas on
   submanifolds and Riemannian products. Then in Section 3, we shall prove the necessity part of the condition, and finally in the last section (Section 4), we shall
   complete the proof of our main theorem by proving the existence and uniqueness.

\section{Preliminaries}

For a natural number $N\geq 2$, denote by $\bbr^N_1$ the Lorentzian space $(\bbr^N,\lagl\cdot,\cdot\ragl_1)$, where the Lorentzian product $\lagl\cdot,\cdot\ragl_1$ is given by
\begin{align}
&\lagl x,y\ragl_1=x^1y^1+\cdots+x^{N-1}y^{N-1}-x^Ny^N,\nonumber\\
&\hs{1.2cm}\forall x=(x^1,\cdots,x^{N-1},x^N),\ y=(y^1,\cdots,y^{N-1},y^N)\in\bbr^N.
\end{align}
Let $S^k$ ($k\geq 1$) be the stand sphere of radius $1$, and $H^{n+p-k}$ the hyperboloid of constant sectional curvature $-1$ which can be isometrically embedded into $\bbr^{n+p-k+1}_1$ as
$$
H^{n+p-k}=\{x=(x^1,\cdots,x^{n+p-k+1})\in\bbr^{n+p-k+1}_1| \lagl x,x\ragl_1=-1,\ x^{n+p-k+1}>0\}.
$$
Therefore we have the standard isometric embedding of the Riemannian product
$$(S^k\times H^{n+p-k},\bar g)\subset \bbr^{k+1}\times\bbr^{n+p-k+1}_1=\bbr^{n+p+2}_1.$$

For each pair of vectors
$$x_1=(x^{1},\cdots,x^{k+1})\in S^k\subset\bbr^{k+1},\quad
x_2=(x^{k+2},\cdots,x^{n+p+2})\in H^{n+p-k}\subset\bbr^{n+p-k+1}_1,$$
Define
$$\xi_{1}=(x_{1},0),\quad
\xi_{2}=(0,x_2).$$
Then $\{\xi_1,\xi_2\}$ is an orthonormal normal frame field of $S^k\times H^{n+p-k}$ in $\bbr^{n+p+2}_1$, and $\xi_1+\xi_2$ is the position vector of $S^k\times H^{n+p-k}$.

At each point $(p_{1},p_{2})$ of $S^k\times H^{n+p-k}$, we have
$$
T_{(p_{1},p_{2})}(S^k\times H^{n+p-k})=T_{p_1}S^k\oplus T_{p_2}H^{n+p-k}.
$$
Thus there is a linear
endomorphism
$$
\bar{\psi}_{(p_{1},p_{2})}:T_{(p_{1},p_{2})}(S^k\times H^{n+p-k})\rightarrow T_{(p_{1},p_{2})}(S^k\times
 H^{n+p-k})
$$
defined by
\be\label{2.1}
\bar{\psi}_{(p_{1},p_{2})}(X_1+X_2)=X_1-X_2,\quad\forall X_1\in T_{p_1}S^k,\ X_2\in T_{p_2}H^{n+p-k}.
\ee
Clearly, $\bar \psi:=\{\bar\psi_{(p_1,p_2)},\ (p_1,p_2)\in S^k\times H^{n+p-k}\}$ defines a differentiable field of bundle endomorphism of $T(S^k\times H^{n+p-k})$, or equivalently, a $(1,1)$-tensor on $S^k\times H^{n+p-k}$. It is easily seen that $\bar \psi$ satisfies the following conditions
\be\label{2.2}
\bar\psi^2=\id,\quad \td D\bar\psi=0,
\ee
where $\td D$ is the induced connection of the Levi-Civita connection $\bar D$ on $S^k\times H^{n+p-k}$. This bundle map $\bar\psi$ is by definition the product structure of $S^k\times H^{n+p-k}$.

\rmk\label{rmk2.1}\rm By defining \be\label{added1}\bar\psi(\xi_1)=\xi_1,\quad \bar\psi(\xi_2)=-\xi_2,\ee we can extend the product structure $\bar\psi$ of $S^k\times H^{n+p-k}$ to be a bundle endmorphism of the trivial bundle
$$T(S^k\times H^{n+p-k})\oplus \bbr\cdot\xi_1\oplus\bbr\cdot\xi_2\equiv (S^k\times H^{n+p-k})\times \bbr^{n+p+2}_1$$
which, composed with the canonical projection, gives exactly the product structure of $\bbr^{k+1}\times\bbr^{n+p-k+1}$. For simplicity, we still use $\bar\psi$ to denote the latter bundle maps. In this case we easily verify that
$$\bar\psi^2=\id,\quad D^0\bar\psi=0,$$
where $D^0$ is the standard flat connection.

 Now let $(M^n,g)$ be a connected Riemannian manifold and $\varphi:M^n\rightarrow S^k\times
  H^{n+p-k}$ an isometric immersion of codimension $p$ with the tangent bundle $TM^n$, the normal bundle $T^\bot_\vfi M^n$ and the second fundamental form $h$. Then we have the following bundle decomposition
  \be\label{2.3}
  T(S^k\times H^{n+p-k})=TS^k\oplus TH^{n+p-k}.
  \ee
  Therefore four bundle maps
  \be\label{2.4}
  f:TM^n
    \rightarrow TM^n, \quad\bar{u}:TM^n\rightarrow T^{\perp}_\vfi M^n,\quad\bar{U}:T^{\perp}_\vfi M^n
    \rightarrow
    TM^n,\quad \bar{\lambda}:T^{\perp}_\vfi M^n\rightarrow T^{\perp}_\vfi M^n
  \ee
are defined by
\be\label{2.5}
\bar{\psi}(\vfi_*X)=\vfi_*f(X)+\bar{u}(X),\quad
\bar{\psi}(\xi)=\vfi_*(\bar{U}(\xi))+\bar{\lambda}(\xi),
\ee
where $X\in TM^n$, $\xi\in\Gamma(T^\bot_\vfi M^n)$.
Use the bundle maps given in \eqref{2.4}, we define the endomorphism
$$\psi:TM^n\oplus T^\bot_\vfi M^n\to TM^n\oplus T^\bot_\vfi M^n$$ by
\be\label{2.5-1}
\psi(X)=f(X)+\bar{u}(X),\quad
\psi(\xi)=\bar{U}(\xi)+\bar{\lambda}(\xi),\quad \forall X\in TM^n,\ \xi\in \xi\in\Gamma(T^\bot_\vfi M^n).
\ee
Then it is not hard to see from \eqref{2.2} that (cf. the proof of
Lemma \ref{l4.4})
\be\label{2.5-2}
\psi^2=\id,\quad \td D\psi=0,
\ee
where $\td D$ is induced by the connection $\bar D$ on $TM^n\oplus T^\bot_\vfi M^n$ given by
\begin{align}
\bar D_XY=&D_XY+h(X,Y),\quad X,Y\in TM^n,\label{2.5-3}\\
\bar D_X\xi=&-A_\xi(X)+D^\bot_X\xi,\quad X\in TM^n,\ \xi\in \Gamma(T^\bot_\vfi M^n),\label{2.5-4}
\end{align}
where, for each $\xi\in \Gamma(E)$, the bundle map $A_\xi:TM^n\to TM^n$
is the Weingarten map of $\vfi$ with respect to $\xi$.

Denote respectively by $D$, $R$, $R^\bot$ and $\bar R$ the Levi-Civita connection of
$(M^n,g)$, the curvature tensor of $D$, the normal curvature of $\vfi$ and the curvature
tensor of the Levi-Civita connection $\bar D$ on $S^k\times H^{n+p-k}$. Then the
fundamental equations of $\vfi$ can be written as follows:
\begin{align}
&\vfi_*(R(X,Y)Z)=(\bar R(X,Y)\vfi_*(Z))^\top +\vfi_*(A_{h(Y,Z)}X)-\vfi_*(A_{h(X,Z)}Y),\label{2.6}\\
&(\bar R(X,Y)\vfi_*(Z))^\bot=(\td D_{X}h)(Y,Z)-(\td D_{Y}h)(X,Z),\label{2.7}\\
&R^{\perp}(X,Y)\xi=(\bar R(X,Y)\xi)^{\perp}+h(A_{\xi}Y,X)-h(A_{\xi}X,Y),\label{2.8}
\end{align}
where $\bar R(X,Y)$ is the curvature operator of the induced bundle $\vfi^*T(S^k\times H^{n+p-k})\to M^n$.

\section{The necessity of the condition in the main theorem}

To prove the main theorem, we begin with the proof of the necessity of the condition.

Let $\vfi:M^n\rightarrow S^k\times H^{n+p-k}$ be as in the last section.

{\lem\label{l2} Denote by $D^0$ the standard flat connection on $\bbr^{n+p+2}_1$. Then

$(1)$ The connection $\bar D$ on the ambient space $S^k\times H^{n+p-k}$ of $\vfi$ is related to the product structure
$\bar\psi$ as follows:
\be\label{2.17}
D^0_{X}Y=\bar{D}_{X}Y-\frac{1}{2}\bar g(X+\bar{\psi}
X,Y)\xi_{1}+\frac{1}{2}\bar g(X-\bar{\psi} X,Y)\xi_{2},\quad X,Y\in T(S^k\times H^{n+p-k}).
\ee

$(2)$ The curvature tensor $\bar R$ of $S^k\times H^{n+p-k}$ is expressed in terms of $\bar\psi$ as
\begin{align}
\bar{R}(X,Y)Z=&\frac{1}{2}(\bar g(\bar{\psi} Y,Z)X+\bar g(Y,Z)\bar{\psi}
X-\bar g(\bar{\psi} X,Z)Y-\bar g(X,Z)\bar{\psi} Y),\nonumber\\
&X,Y,Z\in T(S^k\times H^{n+p-k}).\label{2.18}
\end{align}}

\proof
(1) Denote respectively by $D^{0(1)}$,
$D^{0(2)}$, $\bar{D}^{(1)}$ and $\bar{D}^{(2)}$ the Levi-Civita
connections of $\bbr^{k+1}$, $\bbr^{n+p-k+1}_1$, $S^{k}$ and
$H^{n+p-k}$. Then we have $D^0=D^{0(1)}\times D^{0(2)}$, $\bar D=\bar D^{(1)}\times \bar D^{(2)}$. Let $\bar h^{(1)}$ (resp. $\bar h^{(2)}$) be the
second fundamental form of $S^{k}$ in $\bbr^{k+1}$ (resp. $H^{n+p-k}$
in $\bbr^{n+p-k+1}_1$). Then for all $X=X_{1}+X_{2}$, $Y=Y_{1}+Y_{2}$ with $X_1,Y_{1}\in \Gamma(TS^k)$ and $X_2,Y_{2}\in \Gamma(TH^{n+p-k})$, we have
\begin{align*}
D^0_{X}Y=&D_{X_{1}}^{0(1)}Y_{1}+D_{X_{2}}^{0(2)}Y_{2}\\
=&\bar{D}_{X_{1}}^{(1)}Y_{1}+\bar{D}_{X_{2}}^{(2)}Y_{2}+\bar h^{(1)}(X_{1},Y_{1})+\bar h^{(2)}(X_{2},Y_{2})\\
=&\bar{D}_{X}Y-\bar{g}(X_{1},Y_{1})\xi_{1}+\bar{g}(X_{2},Y_{2})\xi_{2}\\
=&\bar{D}_{X}Y-\frac{1}{2}\bar g(X+\bar{\psi}
X,Y)\xi_{1}+\frac{1}{2}\bar g(X-\bar{\psi} X,Y)\xi_{2},
\end{align*}
where in the last equality we have used the fact that $\bar g(X_1,Y_1)=\bar g(X_1,Y)$, $\bar g(X_2,Y_2)=\bar g(X_2,Y)$ and
$$
X_1=\fr12(X+\bar\psi X), \quad X_2=\fr12(X-\bar\psi X),
$$
which follow easily from the definition of the bundle map $\bar\psi$.

(2) Using \eqref{2.17}, we compute
\begin{align*}
D^{0}_{X}D^{0}_{Y}Z=&D^{0}_{X}(\bar{D}_{Y}Z-\frac{1}{2}\bar
g(Y+\bar{\psi} Y,Z)\xi_{1}+\frac{1}{2}\bar g(Y-\bar{\psi}
Y,Z)\xi_{2})\\
=&\bar{D}_{X}\bar{D}_{Y}Z-\frac{1}{2}\bar g(X+\bar{\psi}
X,\bar{D}_{Y}Z)\xi_{1}+\frac{1}{2}\bar g(X-\bar{\psi}
X,\bar{D}_{Y}Z)\xi_{2})\\
&-\frac{1}{2}\bar
g(\bar{D}_{X}(Y+\bar{\psi}Y),Z)\xi_{1}-\frac{1}{2}\bar
g(Y+\bar{\psi} Y,\bar{D}_{X}Z)\xi_{1}-\frac{1}{2}\bar
g(Y+\bar{\psi} Y,Z)D^0_{X}\xi_{1}\\
 &+\frac{1}{2}\bar
g(\bar{D}_{X}(Y-\bar{\psi}Y),Z)\xi_{2}+\frac{1}{2}\bar
g(Y-\bar{\psi} Y,\bar{D}_{X}Z)\xi_{2}+\frac{1}{2}\bar
g(Y-\bar{\psi} Y,Z)D^0_{X}\xi_{2}.
\end{align*}
We also have
\begin{align*}
D^0_{[X,Y]}Z=\bar D_{[X,Y]}Z-\frac{1}{2}\bar g([X,Y]+\bar{\psi}
[X,Y],Z)\xi_{1}+\frac{1}{2}\bar g([X,Y]-\bar{\psi} [X,Y],Z)\xi_{2}.
\end{align*}
Thus the flatness of $D^0$ and the parallel of $\bar\psi$ show that
\begin{align*}
0=&D^0_XD^0_YZ-D^0_YD^0_XZ-D^0_{[X,Y]}Z\\
=&\bar{R}(X,Y)Z-\frac{1}{2}
\bar g(Y+\bar{\psi} Y,Z)D^0_{X}\xi_{1}+\frac{1}{2}
\bar g(Y-\bar{\psi} Y,Z)D^0_{X}\xi_{2}\\
&\ +\frac{1}{2}
\bar g(X+\bar{\psi} X,Z)D^0_{Y}\xi_{1}-\frac{1}{2}
\bar g(X-\bar{\psi} X,Z)D^0_{Y}\xi_{2}.
\end{align*}
But by the definition of $\xi_1$ and $\xi_2$, we know that
\be\label{3.17-1}
D^0_X\xi_1=X_1=\fr12(X+\bar\psi X), \quad D^0_X\xi_2=X_2=\fr12(X-\bar\psi X).
\ee
Consequently
\begin{align*}
0=&\bar{R}(X,Y)Z-\frac{1}{4}(g(Y+\bar{\psi} Y,Z)(X+\bar{\psi} X)
-\bar g(Y-\bar{\psi} Y,Z)(X-\bar{\psi} X)\\
&\ -g(X+\bar{\psi} X,Z)(Y+\bar{\psi} Y)+g(X-\bar{\psi} X,Z)(Y-\bar{\psi}Y))\\
=&\bar{R}(X,Y)Z-\fr12(\bar g(Y,Z)\bar{\psi} X+\bar
g(\bar{\psi} Y,Z)X-\bar g(X,Z)\bar{\psi} Y-\bar g(\bar{\psi}X,Z)
Y).
\end{align*}
This proves the conclusion (2).
\endproof

{\cor\label{cor} For an isometric immersion $\vfi:M^n\to S^k\times H^{n+p-k}$, if $\bar R(X,Y)$ is the curvature operator of the induced bundle $\vfi^*T(S^k\times H^{n+p-k})$, then the following equations hold:
\begin{align}
\bar{R}(X,Y)\vfi_*(Z)=&\frac{1}{2}(g(f(Y),Z)\vfi_*(X)
+g(Y,Z)\vfi_*f(X)-g(f(X),Z)\vfi_*(Y)-g(X,Z)\vfi_*f(Y))\nonumber\\
&\ +\fr12(g(Y,Z)\bar u(X)-g(X,Z)\bar u(Y)),
\quad X,Y,Z\in \Gamma(TM^n),\label{2.19}\\
\bar{R}(X,Y)\xi=&\frac{1}{2}(\bar g(\bar{\psi} Y,\xi)\vfi_*(X)+\bar g(\vfi_*(Y),\xi)\bar{\psi}
X-\bar g(\bar{\psi} X,\xi)\vfi_*(Y)-\bar g(\vfi_*(X),\xi)\bar{\psi} Y)\nonumber\\
=&\frac{1}{2}\left(\bar g(\bar u(Y),\xi)\vfi_*(X)-\bar g(\bar u(X),\xi)\vfi_*(Y)\right),
\nonumber\\
&\hs{2cm}X,Y\in \Gamma(TM^n),\ \xi\in\Gamma(T^\bot_\vfi M^n). \label{2.20}\end{align}}

The following proposition is a direct consequence of Lemma \ref{l2}, Corollary \ref{cor} together with \eqref{2.6}--\eqref{2.8}:

{\prop\label{p1} The equations of Gauss, Codazzi and Ricci of the immersion $\vfi$ can be rewritten respectively as follows:

$(1)$ Gauss' equation: For all $X,Y,Z\in TM^n$,
\begin{align}
R(X,Y)Z=&A_{h(Y,Z)}X-A_{h(X,Z)}Y\nonumber\\
&\ +\frac{1}{2}(g(f(Y),Z)X+g(Y,Z)f(X)-g(f(X),Z)Y-g(X,Z)f(Y)).\label{2.22}
\end{align}

$(2)$ Codazzi's equation: For all $X,Y,Z\in TM^n$,
\be\label{2.23}
\frac{1}{2}(g(Y,Z)\bar{u}(X)-g(X,Z)\bar{u}(Y))=(\td D_{X}h)(Y,Z)-(\td D_{Y}h)(X,Z).
\ee

$(3)$ Ricci's equation: For all $X,Y\in TM^n$ and $\xi\in T^\bot_\vfi M^n$,
\be\label{2.24}
R^{\bot}(X,Y)\xi=h(A_{\xi}Y,X)-h(A_{\xi}X,Y).
\ee}

Now, by summing up the discussion of this section, the necessity part of the main theorem has been proved and can be formulated as follows:

{\thm\label{t2.1} $(${\bf The necessity}$)$ Let $(M^n,g)$ be a connected Riemannian
manifold of
     dimension $n$, and $\vfi:M^n\to S^k\times H^{n+p-k}$ an isometric immersion with the second fundamental form $h$. Put
     $$
     E=T^\bot_\vfi M^n, \quad \sigma=h,\quad \psi=\bar\psi, \quad \Phi=\id.
     $$
     Then all the equalities \eqref{1.6}--\eqref{1.11} in Theorem \ref{main} hold identically.}

\section{Proof of the existence and the uniqueness}

Having proved the necessity part of the main theorem (Theorem \ref{main}), we shall make use of the techniques of \cite{d1} and \cite{d2} in this section to give a proof of the existence and uniqueness.

Let $N=M^{n}\times R_{1}^{2}$ be the trivial vector bundle with Lorentzian fibre $\bbr^2_1$, and the standard fibre metric is denoted by $g_{-1}$. Let $\{\td\xi_{1},\td\xi_{2}\}$ be an
  orthonormal basis in $\bbr^2_1$ satisfying
\be g_{-1}(\td\xi_{1},\td\xi_{1})=-g_{-1}(\td\xi_{2},\td\xi_{2})=1, \quad g_{-1}({\td\xi_{1},\td\xi_{2}})=0,\ee
 and set $B:=\bar E\oplus N\equiv TM^n\oplus E\oplus N$ as Lorentzian vector bundles with the metric denoted by $\lagl\cdot,\cdot\ragl_B$. Then $\{\td\xi_1,\td\xi_2\}$ can be viewed as an orthonormal frame field of the trivial bundle $N$, thus $\td\xi_1,\td\xi_2$ are sections of $B$.

In terms of the bundle maps $f, u, U, \lambda$ defined in Theorem \ref{main}, we can define an operator $D^B:\Gamma(B)\times \Gamma(TM^n)\to \Gamma(B)$ as follows:
\begin{align}
D^B_{X}Y=&D_{X}Y+\sigma(X,Y)-\frac{1}{2}g(X+f(X),Y)\td\xi_{1}
+\frac{1}{2}g(X-f(X),Y)\td\xi_{2},\quad X,Y\in TM^n,\label{3.3}\\
D^B_{X}\xi=&-A^\sigma_{\xi}X+D^E_{X}\xi-\frac{1}{2}\lagl u(X),\xi\ragl_E(\td\xi_{1}+\td\xi_{2}),\quad X\in TM^n,\quad \xi\in \Gamma(E),\label{3.4}\\
D^B_{X}\td\xi_{1}=&\frac{1}{2}(X+f(X)+u(X)),\quad
D^B_{X}\td\xi_{2}=\frac{1}{2}(X-f(X)-u(X)),\quad X\in TM^n,\label{3.5}
\end{align}
where, for $\xi\in\Gamma(E)$, the bundle endomorphism $A^\sigma_{\xi}$ of $TM^n$ is defined by \eqref{1.10}.

{\lem\label{l3.1} The operator $D^B$ given above is a metric connection on the vector bundle $B$.}

\proof It is a direct verification of the definition of metric connection on vector bundles. \endproof

{\lem\label{l3.2} The bundle maps $f,u,U,\lambda$ satisfy the following conditions
\begin{align}
&g(f(X),Y)=g(f(Y),X), \quad\lagl u(X),\xi\ragl=g(U(\xi),X),\quad \lagl\lambda(\xi),\eta\ragl_E
=\lagl\lambda(\eta),\xi\ragl_E,\label{2.9}\\
&f^{2}+U\circ u =\id _{TM^n}, \quad f\circ U+U\circ \lambda=0,\label{2.10}\\
&u \circ f+\lambda\circ u =0, \quad u\circ U+\lambda^{2}=\id_{E};\label{2.11}\\
  &(\td D_{X}f)(Y)=A^\sigma_{u(Y)}X+U(\sigma(X,Y)),\label{2.12}\\
  &(\td D_Xu)(Y)=\lambda(\sigma (X,Y))-\sigma (X,f(Y));\label{2.13}\\
&(\td D_{X}U)(\xi)=A^\sigma_{\lambda(\xi)}X-f(A^\sigma_{\xi}X),\label{2.14}\\
  &(\td D_{X}\lambda)(\xi)=-\sigma (X,U(\xi))-u(A^\sigma_{\xi}X)\label{2.15},
  \end{align}
where $X,Y\in TM^n$, $\xi,\eta\in \Gamma(E)$.

\proof  We only need to prove the following three statements:

(1) \eqref{2.9} is equivalent to the fact tha $\psi$ is symmetric
respect to the metric $\lagl\cdot,\cdot\ragl_{\bar E}$;

(2) \eqref{2.10} and \eqref{2.11} are equivalent to the fact that $\psi^2=\id$;

(3) \eqref{2.12}--\eqref{2.15} are equivalent to the fact that $\td D\psi=0$.

In fact, the three statements are easily verified via simple computations. For example, the statement (3) is proved as follows:

By \eqref{1.1-1}, \eqref{1.1-2} and the definitions of the $f,u, U, \lambda$, we compute, for any $X,Y\in TM^n$ and $\xi\in\Gamma(E)$,
\begin{align*}
(\tilde{D}_{X}\psi)(Y)=&D^{\bar E}_{X}(\psi(Y))-\psi (D^{\bar E}_{X}Y)\\
=&D^{\bar E}_{X}(f(Y)+u (Y))-\psi ({D}_{X}Y+\sigma (X,Y))\\
=&D_{X}f(Y)-A^\sigma_{u (Y)}X-f({D}_{X}Y)-U(\sigma (X,Y))\\
&\ +\sigma (X,f(Y))+D^E_Xu(Y)-u ({D}_{X}Y)-\lambda(\sigma (X,Y)),
\\
(\tilde{D}_{X}\psi)(\xi)=&D^{\bar E}_{X}(\psi(\xi))-\psi (D^{\bar E}_{X}\xi)\\
=&D^{\bar E}_{X}(U(\xi)+\lambda(\xi))-\psi(-A^\sigma_{\xi}X +\tilde{D}^E_{X}\xi)\\
=&D_{X}U(\xi)-A^\sigma_{\lambda(\xi)}X+fA_{\xi}X-U(D^E_X\xi)\\
&\ +\sigma (X,U(\xi))+D^E_X\lambda(\xi)+u (A^\sigma_{\xi}X)-\lambda(D^E_X\xi).\end{align*}
Note that $\td D\psi=0$ if and only if for all $X,Y\in TM^n$ and $\xi\in \Gamma(E)$,
$$(\tilde{D}_{X}\psi )(Y)=(\tilde{D}_{X}\psi )(\xi)=0.$$
By comparing the tangent and normal components, we easily obtain the
statement (3).\endproof

\lem\label{l4.3} The connection $D^B$ on $B$ defined by \eqref{3.3}--\eqref{3.5} is flat, that is, the curvature operator $R^B$ of $D^B$ vanishes identically.

\proof We need to calculate $R^B(X,Y)\xi$ for arbitrary
field $X, Y\in TM^n$ and $\xi\in\Gamma(B)$. Here we only do the computation for $\xi=Z\in TM^n$, since other cases can be treated similarly. By the definition of the connection $D^B$, we have
\begin{align*}
D^B_{[X,Y]}Z=&D_{[X,Y]}Z+\sigma([X,Y],Z)
-\fr12(g([X,Y],Z)+g(f([X,Y]),Z))\td\xi_1\\
&\ +\fr12(g([X,Y],Z)-g(f([X,Y]),Z))\td\xi_2\\
D^B_{X}D^B_{Y}Z=&D^B_{X}(D_{Y}Z+\sigma(Y,Z)
-\frac{1}{2}g(Y+f(Y),Z)\td\xi_{1}+\frac{1}{2}g(Y-f(Y),Z)\td\xi_{2})\\
=&D_{X}D_{Y}Z-A^\sigma_{\sigma(Y,Z)}X
+D^E_{X}\sigma(Y,Z)+\sigma(X,D_{Y}Z)\\
&\ -\frac{1}{2}(g(X,D_YZ)+g(f(X),D_YZ)+\lagl u(X),\sigma(Y,Z)\ragl_E\\
&\hs{1.5cm} +g(D_XY,Z)+g(D_XfY,Z)+g(Y,D_XZ)+g(f(Y),D_XZ))\td\xi_{1}\\
&\ +\frac{1}{2}(g(X,D_YZ)-g(f(X),D_YZ)-\lagl u(X),\sigma(Y,Z)\ragl_E\\
&\hs{1.5cm} +g(D_XY,Z)-g(D_XfY,Z)+g(Y,D_XZ)-g(f(Y),D_XZ))\td\xi_{2}\\
&\ -\fr12\left(g(Y,Z)f(X)+g(Y,Z)u(X)+g(f(Y),Z)X\right).
\end{align*}
It then follows from \eqref{1.7}, \eqref{1.8} and \eqref{2.12} that
\begin{align*}
R^B(X,Y)Z=&D^B_{X}D^B_{Y}Z-D^B_YD^B_XZ-D^B_{[X,Y]}Z\\
=&R(X,Y)Z+A^\sigma_{\sigma(X,Z)}Y-A^\sigma_{\sigma(Y,Z)}X\\
&\hs{.5cm}+\frac{1}{2}(g(X,Z)f(Y)+g(f(X),Z)Y-g(Y,Z)f(X)-g(f(Y),Z)X)\\
&\ +(\td D_{X}\sigma)(Y,Z)-(\td D_{Y}\sigma)(X,Z)+\frac{1}{2}(g(X,Z)u(Y)-g(Y,Z)u(X))\\
&\ -\fr12(g(A^\sigma_{u(X)}Y-A^\sigma_{u(Y)}X
+(D_{X}f)(Y)-(D_{Y}f)(X),Z))\td\xi_{1}\\
&\ +\fr12(g(A^\sigma_{u(X)}Y-A^\sigma_{u(Y)}X
+(D_{X}f)(Y)-(D_{Y}f)(X),Z))\td\xi_{2}\\
=&0.
\end{align*}
\endproof

Now we extend the bundle map $\psi:\bar E\to \bar E$ to $\td \psi:B\to B$ by defining
\be\label{3.15-1}\td\psi(X)=\psi (X),\text{\ for\ }X\in \bar E;\quad \td\psi(\td\xi_1)=\td\xi_1,\ \td \psi(\td\xi_2)=-\td\xi_2.\ee

{\lem\label{l4.4} The new bundle map $\td\psi$ is parallel with respect to the connection $D^B$, that is, $\td D^B\td\psi=0$, where $\td D^B$ is induced by $D^B$.}

\proof By the definitions of $D^B$ and $\td\psi$, we compute using Lemma \ref{l3.2}
\begin{align*}
(\td D^{B}_{X}\tilde{\psi})(Y)=&D^{B}_{X}(f(Y)+u(Y))-\tilde{\psi}(D^{B}_{X}Y) \\
=&D^{B}_{X}f(Y)+D^{B}_{X}u(Y)\\
&-\tilde{\psi}(D_{X}Y+\sigma(X,Y)-\frac{1}{2}g(X+f(X),Y)\td\xi_{1}+\frac{1}{2}g(X-f(X),Y)\td\xi_{2})\\
=&D_{X}f(Y)+\sigma(X,f(Y))-\frac{1}{2}g(X+f(X),f(Y))\td\xi_{1}+\frac{1}{2}g(X-f(X),f(Y))\td\xi_{2}\\
&-A^\sigma_{u(Y)}X+D^E_{X}u(Y)-\frac{1}{2}\lagl u(X),u(Y)\ragl_E(\td\xi_{1}+\td\xi_{2})\\
&-f(D_{X}Y)-u(D_{X}Y)-U(\sigma(X,Y))-\lambda(\sigma(X,Y))\\
&+\frac{1}{2}g(X+f(X),f(Y))\td\xi_{1}+\frac{1}{2}g(X-f(X),f(Y))\td\xi_{2}\\
=&(\td D_{X}f)(Y)-A^\sigma_{u(Y)}X-U(\sigma(X,Y))\\
&+(\tilde{D}_{X}u)(Y)+\sigma(X,f(Y))-\lambda(\sigma(X,Y))\\
&+\frac{1}{2}(g(-f(X)-f^{2}X-U(u(X))+X+f(X),Y)\td\xi_{1}\\
&+\frac{1}{2}(g(f(X)-f^{2}X-U(u(X))+X-f(X),Y)\td\xi_{2} =0;\\
(\td D^{B}_{X}\tilde{\psi})(\xi)=&D^{B}_{X}\tilde{\psi}(\xi)-\tilde{\psi}(D^{B}_{X}\xi)\\
=&D^{B}_{X}(U(\xi)+\lambda(\xi))-\tilde{\psi}(-A^\sigma_{\xi}X +D^E_{X}\xi-\frac{1}{2}\lagl u(X),\xi\ragl_E(\td\xi_{1}+\td\xi_{2}))\\
=&D_{X}U(\xi)+\sigma(X,U(\xi))-\frac{1}{2}g(X+f(X),U(\xi))\td\xi_{1}+\frac{1}{2}g(X-f(X),U(\xi))\td\xi_{2}\\
&-A^\sigma_{\lambda(\xi)}X+D^E_{X}\lambda(\xi)-\frac{1}{2}\lagl u(X),\lambda(\xi)\ragl_E(\td\xi_{1}+\td\xi_{2})\\
&+f(A^\sigma_{\xi}X)+u(A^\sigma_{\xi}X)-U(D^E_{X}\xi)-\lambda(D^E_{X}\xi)\\
&+\frac{1}{2}\lagl u(X),\xi\ragl_E\td\xi_{1}-\frac{1}{2}\lagl u(X),\xi\ragl_E\td\xi_{2}\\
=&(\td D_{X}U)(\xi)-A^\sigma_{\lambda(\xi)}X+f(A^\sigma_{\xi}X)\\
 &+(\tilde{D}_{X}\lambda)(\xi)+\sigma(X,U(\xi))+u(A^\sigma_{\xi}X)\\
&-\frac{1}{2}\lagl u(f(X))+\lambda(u(X)),\xi\ragl_E\td\xi_{1}\\
&-\frac{1}{2}\lagl u(f(X))+\lambda(u(X)),\xi\ragl_E\td\xi_{2}=0;\\
(\td D^{B}_{X}\tilde{\psi})(\td\xi_{1}) =&D^{B}_{X}\td\psi(\td\xi_{1})-\tilde{\psi}(D^{B}_{X}\td\xi_{1})
=D^{B}_{X}\td\xi_{1}-\fr12\td\psi(X+f(X)+u(X))\\
=&\frac{1}{2}(X+f(X)+u(X))-\fr12(f(X)+f^2(X)+u(X)+u(f(X))+U(u(X)))+\lambda(u(X)))\\
=&\frac{1}{2}(X-f^{2}X-U(u(X))-u(f(X))-\lambda(u(X)))=0;\\
(\td D^{B}_{X}\tilde{\psi})(\td\xi_{2}) =&D^{B}_{X}\td\psi(\td\xi_{2})-\tilde{\psi}(D^{B}_{X}\td\xi_{2})
=-D^B_X\td\xi_2-\fr12\td\psi(X-f(X)-u(X))\\
=&-\frac{1}{2}(X-f(X)-u(X))-\fr12(f(X)-f^2X+u(X)-u(f(X))-U(u(X))-\lambda(u(X)))\\
=&\frac{1}{2}(-X+f^{2}X+U(u(X))+u(f(X))+\lambda(u(X)))=0.
\end{align*}
\endproof

Let $B_{1}$ and $B_2$ be the subbundle $B$ defined respectively by
\be\label{3.13}B_{1}=\{b\in B;\ \tilde{\psi}b=b\},\quad B_2=\{b\in
B;\ \tilde{\psi}b=-b\}.\ee
Clearly, $B_1$ is Riemannian while $B_2$ is Lorentzian and, by \eqref{3.15-1}, the two line bundles $\bbr\td\xi_1$ and $\bbr\td\xi_2$ are respectively the subbundles of $B_1$ and $B_2$.
Since $\td\psi$ is symmetric and $\td\psi^2=\id$, $\psi\neq\pm\id$, it is easily seen that $B_1$ and $B_2$ are orthogonal each other and, if $k+1$ (resp. $n+p-k+1$) is the rank of $B_1$ (resp. $B_2$), then $1\leq k\leq n+p-1$ (resp. $1\leq n+p-k\leq n+p-1$).

{\lem\label{l4.5} Under the assumption of Theorem \ref{main}, the subbundle $B_1$ $($resp. $B_2)$ has a parallel and orthonormal frame field $\{s_1,\cdots,s_{k+1}\}$ $($resp. $\{s_{k+2},\cdots,s_{n+p+2}\})$ defined on $M^n$, where $s_{n+p+2}$ is a time-like unit section, that is, $\bar g(s_{n+p+2},s_{n+p+2})=-1$. Furthermore, $$\{s_1,\cdots,s_{k+1},s_{k+2},\cdots,s_{n+p+2}\}$$ is a parallel $($Lorentzian$)$ orthonormal frame field of the bundle $B$ globally defined on $M^n$.}

\proof To start, we fix one point $p_0\in M^n$ and choose, at $p_0$, one orthonormal basis $\{b_1,\cdots,b_{k+1}\}$ of $B_1$ and one Lorentzian orthonormal basis $\{b_{n+p-k+2}\}$ of $B_2$, respectively, such that $b_{k+1}=\td\xi_1$, and $b_{n+p+2}=\td\xi_2$. Since $M^n$ is simply connected and $B$ is flat by Lemma \ref{l4.3}, there exists one globally defined parallel frame field $s:=\{s_1,\cdots,s_{k+1},s_{k+2},\cdots,s_{n+p+2}\}$ of $B$, such that $s_i(p_0)=b_i$, $1\leq i\leq n+p+2$.

On the other hand, since $\td\psi$ is parallel by Lemma \ref{l4.4}, it easily follows that $s_1,\cdots,s_{k+1}$ (resp. $s_{k+2},\cdots,s_{n+p+2}$) are parallel sections of the subbundle $B_1$ (resp. $B_2$). This
completes the proof of Lemma \ref{l4.5}.\endproof

Now we are ready to prove the existence and the uniqueness separately as follows:

{\bf 1. Proof of the existence}

Step 1. The basic formulation of a map $\vfi:M^n\to S^k\times H^{n+p-k}$.

Denote by $B_p$ the fibre of $B$ at each point $p\in M^n$, and use similar notations for other bundles. Let $s=\{s_1,\cdots,s_{k+1},s_{k+2},\cdots,s_{n+p+2}\}$ be the parallel frame field of $B$ given by Lemma \ref{l4.5}. Then we have a natural isomorphism $\Psi:B\to M^n\times \bbr^{n+p+2}_1$ of vector bundles defined by
\be\label{4.14}
\Psi(b)=(p,(b^1,\cdots,b^{n+p+2})),\quad \forall\, b=\sum_{i=1}^{n+p+2}b^is_i\in B_p,\ p\in M^n.
\ee
Since the frame field $s$ is orthonormal with $s_{n+p+2}$ time-like, it follows that $\Psi$ is metric-preserving and
\be\label{4.15} b^i=\lagl b,s_i\ragl_B,\ \text{for\ }i=1,\cdots,n+p+1; \quad b^{n+p+2}=-\lagl b,s_{n+p+2}\ragl_B.\ee
Note that the frame field $s$ is parallel, so $\Psi$ is also a bundle map preserving the connections.

Define $\vfi:=\pi_2\circ\Psi(\td\xi_1+\td\xi_2)$, where $\pi_2:M^n\times\bbr^{n+p+2}_1\to \bbr^{n+p+2}_1$ is the projection onto the second factor. Then it easily seen that $\vfi(M^n)\subset S^k\times H^{n+p-k}$. Consequently, $\vfi$ is a map from $M^n$ into the product space $S^k\times H^{n+p-k}$. Write
\be\label{4.16} \td\xi_1=\sum_{i=1}^{k+1}x^is_i,\quad \td\xi_2=\sum_{j=k+2}^{n+p+2}y^js_j,
\ee
where $x^i,y^j\in C^\infty(M^n)$, $1\leq i\leq k+1$, $k+2\leq j\leq n+p+2$.
Then, as an $\bbr^{n+p+2}$-valued function,
\be\label{4.17}\vfi=(x^1,\cdots,x^{k+1},y^{k+2},\cdots,y^{n+p+2}).\ee

Step 2. $\vfi:M^n\to S^k\times H^{n+p-k}$ is an isometric immersion.

Here we only need to prove that $\vfi$ is isometric, that is, $\bar g(\vfi_*(X),\vfi_*(Y))=g(X,Y)$ for all $X,Y\in TM^n$.
Since $s$ is parallel and $D^B$ is compatible with the metric $\lagl\cdot,\cdot\ragl_B$, it is seen from \eqref{4.15} and \eqref{4.17} that,
for $X\in TM^n$,
\begin{align}
\vfi_*(X)=&X(\vfi)=(X(x^1),\cdots,X(x^{k+1}),X(y^{k+2}),\cdots,X(y^{n+p+2}))\nonumber\\
=&(X\lagl\td\xi_1,s_1\ragl_B,\cdots, X\lagl\td\xi_1,s_{k+1}\ragl_B,X\lagl\td\xi_2,s_{k+2}\ragl_B,\cdots,-X\bar \lagl\td\xi_2,s_{n+p+2}\ragl_B)\nonumber\\
=&(\lagl D^B_X\td\xi_1,s_1\ragl_B,\cdots,
\lagl D^B_X\td\xi_1,s_{k+1}\ragl_B,\lagl D^B_X\td\xi_2,s_{k+2}\ragl_B,\cdots,
-\lagl D^B_X\td\xi_2,s_{n+p+2}\ragl_B)\nonumber\\
=&\Psi(D^B_X\td\xi_1+D^B_X\td\xi_2).\label{4.19}
\end{align}
On the other hand, from \eqref{3.5}, it follows that
$$
D^B_X\td\xi_1+D^B_X\td\xi_2=\fr12(X+f(X)+u(X))+\fr12(X-f(X)-u(X))=X.
$$
Consequently, $\vfi_*$ can be identified with $\Psi|_{TM^n}$, that is
\be\label{4.19-1}
\vfi_*(X)=\Psi(X),\quad\text{for all }X\in TM^n.
\ee

Thus by the fact that $\Psi$ is metric-preserving, we obtain for $X,Y\in TM^n$
$$
\bar g(\vfi_*(X),\vfi_*(Y))=\lagl\vfi_*(X),\vfi_*(Y)\ragl_1 =\lagl\Psi(X),\Psi(Y)\ragl_1=\lagl X,Y\ragl_B=g(X,Y).
$$

Step 3. The normal bundle $T^\bot_\vfi M^n=\Psi(E)$.

For an arbitrary vector $\xi\in E_p$ ($p\in M^n$), we have by \eqref{4.14} and \eqref{4.15}
$$
\Psi(\xi)=(p,(\lagl\xi, s_{1}\ragl_B, \cdots, \lagl\xi,s_{k+1}\ragl_B,\lagl\xi, s_{k+2}\ragl_B,  \cdots,
\lagl\xi,s_{n+p+1}\ragl_B,-\lagl\xi,s_{n+p+2}\ragl_B)).
$$
Thus for each $X\in T_pM^n$, we use \eqref{4.19-1} to get
$$\langle\varphi_{*}X,\Psi(\xi)\rangle_1=\langle\Psi(X),\Psi(\xi)\rangle_1=\lagl X,\xi\ragl_B=0.$$
Hence $\Psi(E)\subset T^\bot_\vfi M^n$. Since $\Psi$ is bundle-isomorphic and $\rank E=\rank T^\bot_\vfi M^n=p$,
we obtain that $\Psi(E)=T^\bot_\vfi M^n$.

Step 4. $\vfi:M^n\to S^k\times H^{n+p-k}$ satisfies all the conditions of the main theorem.

Let $h$ be the second fundamental form of $\vfi$ and denote $\Phi=\Psi|_E$. Then $\Phi:E\to T^\bot_\vfi M^n$ is clearly an isometric bundle map. Since $\Psi$ is connection-preserving, we obtain by \eqref{4.19-1} that,
for $X,Y\in\Gamma(TM^n)$ and $\xi\in T^\bot_\vfi M^n$,
\begin{align*}\bar g(h(X,Y),\xi)=&\bar g(D^0_X\vfi_*(Y),\xi)=\bar g(D^0_X\Psi(Y),\xi)
=\bar g(\Psi(D^B_XY),\xi)=\lagl D^B_XY,\Psi^{-1}(\xi)\ragl_B\\
=&\lagl D_XY+\sigma(X,Y)-\fr12 g(X+f(X),Y)\td\xi_1+\fr12 g(X-f(X),Y)\td\xi_2,\Psi^{-1}(\xi)\ragl_B\\
=&\lagl\sigma(X,Y),\Psi^{-1}(\xi)\ragl_B=\bar g(\Psi(\sigma(X,Y)),\xi)=\bar g(\Phi(\sigma(X,Y)),\xi).
\end{align*}
Then, by the arbitrariness of $\xi\in T^\bot_\vfi M^n$ and $X,Y\in TM^n$, it follows that $h=\Phi\circ\sigma$.

Now by the definition of the isometric bundle map $\Psi:B\to M^n\times\bbr^{n+p+2}_1$, we clearly see that the extended bundle maps $\bar\psi$ and $\td\psi$ given respectively in Remark \ref{rmk2.1} and \eqref{3.15-1} satisfy $\bar\psi\circ\Psi=\Psi\circ\psi$. Therefore, for any $X\in TM^n$ and $\xi\in E$, we have by \eqref{4.19} and the parallel of $\psi$
\begin{align*}
\bar{\psi}(\vfi_*(X))=&\bar\psi(\Psi(D^B_X\td\xi_1+D^B_X\td\xi_2))
=\Psi(\td\psi(D^B_X\td\xi_1+D^B_X\td\xi_2))\\
=&\Psi(D^B_X(\td\psi(\td\xi_1))+D^B_X(\td\psi(\td\xi_2))) =\Psi(D^B_X\td\xi_1-D^B_X\td\xi_2)\\
=&\fr12\Psi(X+f(X)+u(X)-X+f(X)+u(X))\\
=&\Psi(f(X))+\Psi(u(X))
=\vfi_*(f(X))+\Phi(u(X)),\\
\bar{\psi}(\Phi(\xi))=&\bar{\psi}(\Psi(\xi))=\Psi(\td\psi(\xi))=\Psi(U(\xi)+\lambda(\xi))\\
=&\vfi_*(U(\xi))+\Phi(\lambda(\xi)).
\end{align*}

{\bf 2. Proof of the uniqueness}

Finally we prove that the isometric immersion is unique up to
isometries of $S^k\times H^{n+p-k}$. Let $\bar\varphi:M^{n}\rightarrow S^k\times H^{n+p-k}$ be another isometric
immersion with isometric bundle map $\bar\Phi:E\to T^\bot_{\bar\vfi} M^n$ of the given metric bundle $E\to M^n$ onto the normal bundle $T^\bot_{\bar\vfi}M^n\to M^n$ of $\bar\vfi$, satisfying all the conditions in Theorem \ref{main}. In particular, $\bar h:=\bar\Phi\circ\sigma$ is the second fundamental form of $\bar\vfi$.

In terms of \eqref{4.17}, we write
$$
\xi_1=(x^1,\cdots,x^{k+1},0,\cdots,0),\quad \xi_2=(0,\cdots,0,y^{k+2},\cdots,y^{n+p+2}).
$$
Similarly, if
$$\bar\vfi=(\bar x^1,\cdots,\bar x^{k+1},\bar y^{k+2},\cdots,\bar y^{n+p+2}),$$
we write
$$
\bar\xi_1=(\bar x^1,\cdots,\bar x^{k+1},0,\cdots,0),\quad \bar\xi_2=(0,\cdots,0,\bar y^{k+2},\cdots,\bar y^{n+p+2}).
$$

Similar to the bundle map $\Psi:B\to M^n\times\bbr^{n+p+2}_1$, we define another bundle map $\bar\Psi:B\to M^n\times\bbr^{n+p+2}_1$ by
\be
\bar\Psi|_{TM^n}=\bar\vfi_*,\quad \bar\Psi|_E=\bar\Phi,\quad \bar\Psi(\td\xi_1)=\bar\xi_1,\ \bar\Psi(\td\xi_2)=\bar\xi_2.
\ee
Then we have
{\lem\label{l4.6} The bundle map $\bar\Psi:B\to M^n\times\bbr^{n+p+2}_1$ defined above keeps invariant both the metrics and the connections.}

\proof That $\bar\Psi$ keeps the fibre metrics invariant follows directly from the fact that
$$\bar\vfi_*:TM^n\to \bar\vfi_*(TM^n),\quad \bar\Phi:E\to T^\bot_{\bar\vfi}M^n$$
are metric-preserving and that both $\{\td\xi_1,\td\xi_2\}$ and $\{\bar\xi_1,\bar\xi_2\}$ are orthonormal frame fields.

To prove that $\bar\Psi$ keeps the connections invariant we first note by the Gauss Formula that
$$
\bar D_X\bar\vfi_*(Y)=\bar\vfi_*(D_XY)+\bar h(X,Y)
=\bar\Psi(D_XY)+\bar \Psi(\sigma(X,Y))
=\bar\Psi(D_XY+\sigma(X,Y)).
$$
Hence from \eqref{2.17} and \eqref{3.3}, it follows that, for all $X,Y\in TM^n$
\begin{align*}
D^0_X\bar\Psi(Y)=&D^0_X\vfi_*(Y)=\bar D_X\vfi_*(Y)-\fr12(g(X+f(X),Y))\bar\xi_1
+\fr12(g(X-f(X),Y))\bar\xi_2\\
=&\bar\Psi(D_XY+\sigma(X,Y))
-\fr12(g(X+f(X),Y))\bar\Psi(\td\xi_1)
+\fr12(g(X-f(X),Y))\bar\Psi(\td\xi_2)\\
=&\bar\Psi(D_XY+\sigma(X,Y)-\fr12\lagl X+f(X),Y\ragl_B\td\xi_1
+\fr12\lagl X-f(X),Y\ragl_B\td\xi_2)\\
=&\bar\Psi(D^B_XY)
\end{align*}
which gives $D^0\circ\bar\Psi=\bar\Psi\circ D^B$ on $TM^n$.

Let $\bar A_{\bar\xi}:TM^n\to TM^n$ be the Weingarten map of $\bar\vfi$ and $\{e_1,\cdots,e_n\}$ a fixed orthonormal frame field of $TM^n$. Then by \eqref{1.10} we have for any $X\in TM^n,\xi\in E$,
\begin{align*}
\bar\vfi_*(\bar A_{\bar\Psi(\xi)}X)=&\bar\vfi_*\left(\sum_i g(\bar A_{\bar\Psi(\xi)}X,e_i)e_i\right)
=\bar\vfi_*\left(\sum_i\bar g(\bar h(X,e_i),\bar\Psi(\xi))e_i\right)\\
=&\bar\vfi_*\left(\sum_i\bar g(\bar\Psi(\sigma(X,e_i)),\bar\Psi(\xi))e_i\right)
=\bar\vfi_*\left(\sum_i\lagl\sigma(X,e_i),\xi\ragl_E e_i\right)
=\bar\Psi\left(\sum_i g(A^\sigma_\xi X,e_i)e_i\right)\\
=&\bar\Psi(A^\sigma_\xi X).
\end{align*}
From \eqref{1.11}, \eqref{2.17}, \eqref{3.4}, together with the Weingarten formula of $\bar\vfi$
and the fact that the bundle map $\bar\Phi:E\to T^\bot_{\bar\vfi}M^n$ keeps the connections invariant, it follows that
\begin{align*}
D^0_X\bar\Psi(\xi)=&D^0_X\bar\Phi(\xi)=\bar D_X\bar\Phi(\xi)-\fr12\bar g(\bar \psi(\vfi_*(X)),\bar\Phi(\xi))(\bar\xi_1
+\bar\xi_2)\\
=&-\bar\vfi_*(\bar A_{\bar\Psi(\xi)}X)+\bar D^\bot_X\bar\Phi(\xi)
-\fr12\bar g(\vfi_*(f(X))+\bar\Phi(u(X)),\bar\Phi(\xi))
(\bar\xi_1
+\bar\xi_2)\\
=&-\sum_ig(\bar A_{\bar\Psi(\xi)}X,e_i)\vfi_*(e_i)+\bar\Phi(D^E_X(\xi))
-\fr12\bar\Psi(\bar g(\bar\Phi(u(X)),\bar\Phi(\xi))(\td\xi_1
+\td\xi_2))\\
=&-\sum_i\bar g(\bar h(X,e_i),\bar\Psi(\xi))\bar\Psi(e_i)+\bar\Phi(D^E_X(\xi))
-\fr12\bar\Psi(\lagl u(X),\xi\ragl_E(\td\xi_1
+\td\xi_2))\\
=&-\sum_i\bar g(\bar\Psi(\sigma(X,e_i),\bar\Psi(\xi)))\bar\Psi(e_i)
+\bar\Phi(D^E_X(\xi))
-\fr12\bar\Psi(\lagl u(X),\xi\ragl_E(\td\xi_1
+\td\xi_2))\\
=&\bar\Psi\left(-A^\sigma_\xi(X)+D^E_X(\xi) -\fr12\lagl u(X),\xi\ragl_E(\td\xi_1+\td\xi_2)\right)\\
=&\bar\Psi(D^B_X\xi)
\end{align*}
showing that the equality $D^0\circ\bar\Psi=\bar\Psi\circ D^B$ also holds on $E$.

Furthermore, by \eqref{3.17-1} and \eqref{3.5}, for all $X\in TM^n$
$$D^0_X\bar\Psi(\td\xi_1)=D^0_X\bar\xi_1
=\fr12(\vfi_*(X+f(X))+\bar\Phi(u(X)))
=\fr12(\bar\Psi(X+f(X)+u(X)))=\fr12(\bar\Psi(D^B_X\td\xi_1)).$$
In the same way, $D^0_X\bar\Psi(\td\xi_2)=\bar\Psi(D^B_X\td\xi_2)$.

Now by summing up we conclude that $D^0\circ\bar\Psi=\bar\Psi\circ D^B$ identically on $B$.\endproof

Set $T=\bar\Psi\circ\Psi^{-1}$. Then it is clear that $T$ uniquely defines a smooth map $T:M^n\to O(n+p+2,1)$ from $M^n$ into the Lorentzian orthogonal group $O(n+p+2,1)$, which is viewed as a family of isometries of $\bbr^{n+p+2}_1$ parameterized by $M^n$ satisfying $\bar\vfi=T(\vfi)$.

On the other hand, we know by \eqref{1.11} that $\Psi\circ\td\psi=\bar\psi\circ\Psi$,
$\bar\Psi\circ\td\psi=\bar\psi\circ\bar\Psi$. It follows that $T\circ\bar\psi=\bar\psi\circ T$ and thus, for each $p\in M^n$,
the value $T(p)$ of the map $T:M^n\to O(n+p+2,1)$ is an isometry of the product space $S^k\times
H^{n+p-k}$. Finally, since both $\Psi$ and $\bar\Psi$ are parallel as bundle maps, we obtain by the definition of the standard connection $D^0$ on the trivial bundle $M^n\times\bbr^{n+p+2}_1$ that $dT\equiv 0$. Therefore
$T$ is constant thus a fixed isometry of $S^k\times
H^{n+p-k}$, completing the proof of the uniqueness.


\ed
Define
\be\label{4.20}\xi_1:=\Psi(\td\xi_1)=(x^1,\cdots,x^{k+1},0,\cdots,0),\quad
\xi_2:=\Psi(\td\xi_1)=(0,\cdots,0,y^{k+2},\cdots,y^{n+p+2}).
\ee